%% file: sucag_cdc_fin.tex
\newif\ifconfver
\confverfalse     

\newif\ifplainver  
\plainvertrue

\ifplainver
    \confverfalse   
\fi

\ifconfver

\documentclass[conference]{IEEEtran}

\else
    \ifplainver
        \documentclass[11pt]{article}
        \usepackage{fullpage}
    \else
        \documentclass[11pt,draftcls,onecolumn]{IEEEtran}
    \fi
\fi
\usepackage{calc,amsfonts,amssymb,amsmath,bm,url,color,theorem,graphicx,cite,shortcuts_FW}
\usepackage{psfrag,subfigure,float,hyperref}
\usepackage{algorithm}
\usepackage{algorithmic}
\usepackage{multirow, enumitem}

\usepackage[font=small,labelfont=bf,tableposition=top]{caption}

\definecolor{orange}{RGB}{255,107,0}

\newtheorem{Lemma}{Lemma}

\newtheorem{Theorem}{Theorem}
\newtheorem{Assumption}{A\hspace{-.1cm}}

\usepackage{pgfplots}
\pgfplotsset{compat=1.3}
\usepackage{tikz}
\usetikzlibrary{arrows,shapes,calc,tikzmark,backgrounds,matrix,decorations.markings}
\usepgfplotslibrary{fillbetween}

\usetikzlibrary{shadows.blur}
\definecolor{asuorange}{rgb}{1,0.699,0.0625}
\definecolor{asured}{rgb}{0.598,0,0.199}
\definecolor{asuborder}{rgb}{0.953,0.484,0}
\definecolor{asugrey}{rgb}{0.309,0.332,0.340}
\definecolor{asublue}{rgb}{0,0.555,0.836}
\definecolor{asugold}{rgb}{1,0.777,0.008}

\definecolor{lavander}{cmyk}{0,0.48,0,0}
\definecolor{violet}{cmyk}{0.79,0.88,0,0}
\definecolor{burntorange}{cmyk}{0,0.52,1,0}

\tikzstyle{stubborn}=[draw,circle, black!80, fill=black!40,
                       text=violet,minimum width=10pt]
\tikzstyle{superpeers}=[draw,circle, asublue!80!white, fill = asublue!50!white,
                       text=violet,minimum width=10pt]
\tikzstyle{susceptible}=[draw,circle, left color = red, color = red,
                       text=violet,minimum width=10pt]
\tikzstyle{perturb}=[draw,circle,burntorange, left color=blue,
                       text=violet,minimum width=10pt]
\tikzstyle{legend_general}=[rectangle, rounded corners, thin,
                           black, fill= white, draw, text=black,
                           minimum width=2.5cm, minimum height=0.8cm]
\tikzstyle{legend_graph}=[rectangle, rounded corners, thin,
                           black, fill= white, draw, text=black,
                           minimum width=1cm, minimum height=0.5cm]
\tikzstyle{legend_fw}=[rectangle, rounded corners, very thin,
                           black, fill= asuorange!8, draw, text=black]

\usepackage{etoolbox}
\makeatletter
\patchcmd{\@makecaption}
  {\scshape}
  {}
  {}
  {}
\makeatletter
\patchcmd{\@makecaption}
  {\\}
  {.\ }
  {}
  {}
\makeatother

\newcommand{\hide}[1]{}

\definecolor{mygreen}{rgb}{0.0, 0.5, 0.0}

\def\nfc#1{{\color{red} NF: #1}}

\begin{document}

\bibliographystyle{IEEEtran}

\newcommand{\papertitle}
{\vspace{.48cm}
SUCAG: Stochastic Unbiased Curvature-aided Gradient Method 
for Distributed Optimization
}

\newcommand{\paperabstract}{We propose and analyze a new stochastic gradient method, which we call \emph{Stochastic Unbiased Curvature-aided Gradient} (SUCAG),
for finite sum optimization problems. 
SUCAG constitutes an \emph{unbiased} total gradient tracking technique  
that uses Hessian information to accelerate convergence. 
We analyze our method under the general asynchronous model of computation, in which each function is selected infinitely often with 
possibly unbounded (but sublinear) delay. 
For strongly convex problems, we establish linear convergence
for the SUCAG method. When the initialization point is sufficiently close
to the optimal solution, the established convergence rate is 
only dependent on the condition number of the problem,
making it strictly faster than the known rate for the SAGA method.

%

Furthermore, we describe a Markov-driven approach of implementing 
the SUCAG method
in a distributed asynchronous multi-agent setting, via 
gossiping along a random walk
on an undirected communication graph. We show that our analysis applies as long as the graph is connected and, notably, establishes an asymptotic linear convergence rate that \emph{is robust to the graph topology}. 
Numerical results demonstrate
the merits of our algorithm over existing methods.}


\ifplainver


    \title{\papertitle}

    \author{
    Hoi-To Wai, Nikolaos M. Freris, Angelia Nedi\'c and Anna Scaglione
    \thanks{This work is supported by the NSF, under grants NSF CCF-BSF 1714672, CCF-1717207 and CCF-1717391. H.-T.~Wai, A.~Nedi\'c, A.~Scaglione are with the School of Electrical, Computer and Energy Engineering, Arizona State University, Tempe, AZ 85281, USA. E-mails: \texttt{\{htwai,Angelia.Nedich,Anna.Scaglione\}@asu.edu}. N.~M.~Freris is with Division of Engineering, New York University Abu Dhabi \& NYU Tandon School of Engineering. E-mail: \texttt{nf47}@nyu.edu.}
    }

    \maketitle

\else
    \title{\papertitle}

    \ifconfver \else {\linespread{1.1} \rm \fi

    \author{\IEEEauthorblockN{Hoi-To Wai\IEEEauthorrefmark{1}, Nikolaos M. Freris\IEEEauthorrefmark{2}, Angelia Nedi\'c\IEEEauthorrefmark{1} and Anna Scaglione\IEEEauthorrefmark{1}}
\thanks{This work was supported by the National Science Foundation (NSF), under grants NSF CCF-BSF 1714672 and CCF-1717207.}
\IEEEauthorblockA{\IEEEauthorrefmark{1}School of Electrical, Computer and Energy Engineering, Arizona State University, Tempe, AZ 85281, USA.}
\IEEEauthorblockA{\IEEEauthorrefmark{2}Division of Engineering, New York University Abu Dhabi \& NYU Tandon School of Engineering\\
Emails: \texttt{\{htwai,angelia.nedich,Anna.Scaglione\}@asu.edu}, \texttt{nf47@nyu.edu}}
         }

    \maketitle

    \ifconfver \else
        \begin{center} \vspace*{-2\baselineskip}
        \end{center}
    \fi

    \ifconfver \else \IEEEpeerreviewmaketitle} \fi

 \fi
 
     \begin{abstract}
\paperabstract\end{abstract}

\ifplainver
\else
    \begin{IEEEkeywords}\vspace{-0.0cm}
       Distributed optimization, Incremental methods, Asynchronous algorithms, Randomized algorithms, Multi-agent systems, Machine learning.
    \end{IEEEkeywords}
\fi

\ifconfver \else
    \ifplainver \else
        \newpage
\fi \fi

\section{Introduction}
We consider the finite sum optimization problem:
\beq \label{eq:opt}
\min_{ \prm \in \RR^d }~ F(\prm) 
~~\text{where}~~F(\prm) \eqdef \frac{1}{N} \sum_{i=1}^N f_i (\prm)  \eqs.
\eeq
Each component function 
$f_i : \RR^d \rightarrow \RR$ is assumed to be convex and twice continuously differentiable with Lipschitz continuous gradients, 
while the sum function $F(\prm)$ is assumed to be strongly convex. 
Such problem customarily arises in several multi-agent systems applications \cite{nedic2009distributed}, e.g., communication networks~\cite{lin2006tutorial}, sensor networks~\cite{freris2010fundamentals}, and distributed learning~\cite{boyd2011distributed}. 
In a distributed setting, $f_i (\cdot)$ is the \emph{private} 
function known only to agent $i$; for example, in a learning task such as empirical risk minimization (ERM)~\cite{vapnik1999overview}, this corresponds to a subset of data gathered by agent $i$. 

This paper 
develops an efficient distributed algorithm for \eqref{eq:opt} under a general communication topology. In particular, we adopt 
a stochastic asynchronous setting that captures activation of agents in networked decision and control applications. 
To this end, a large body of prior work has focused on developing distributed algorithms based on the average consensus protocol \cite{Tsitsiklis1984}, e.g., 
as introduced by {Nedi\'{c}} \emph{et al.} \cite{nedic2001distributed}.
However, the convergence rate for this class of methods is limited by 
the diameter of the communication graph \cite{scaman2017optimal}.  
To develop an algorithm that is robust to the graph structure, 
we consider here an alternative pathway that is 
akin to that described in \cite{johansson2009randomized,ram2009incremental}, 
for utilizing stochastic algorithms for finite sum minimization 
via pairwise exchange of information on a graph \cite{nedic2009distributed,RK}. 
Inspired by the classical work of Robbins and Monro~\cite{robbins1951stochastic}, stochastic algorithms have been actively pursued for machine learning operations over the last decade, with notable recent examples 
enlisting stochastic average gradient (SAG)\cite{SAG} and its unbiased 
variant (SAGA)~\cite{SAGA}, and stochastic variance reduction gradient (SVRG)~\cite{SVRG}. All these methods (including the one we propose here)
are \emph{incremental}~\cite{nedic2001incremental} and operate using the following key steps, at each iteration:   
a) a \emph{random} component function is selected, and its gradient is evaluated at the current iterate, b) the obtained gradient value is aggregated into the memory, and c) a new iterate is obtained using the aggregated gradient. Additionally, they require storing the gradient values to allow for an asynchronous implementation. 

Inspired by \cite{johansson2009randomized,ram2009incremental}, 
we make an observation that the \emph{random selection-aggregate-update} paradigm can be realized
in a fully distributed fashion via a random walk on a graph, where a selected agent forwards information to one of its neighbors, \ie  essentially implementing a Markov-driven token passing protocol;
cf.~Section~\ref{sec:dist}. 
Most notably, while 
stochastic algorithms such as SAG and SAGA are known 
to 
converge linearly for strongly convex problems~\cite{SAG,SAGA}, 
their rates depend on the number of agents $N$.  
For example, SAGA requires 
${\cal O}( ( N + \kappa(F) ) \log (1/\epsilon) )$ steps in expectation  
to reach an $\epsilon$-optimal solution, where $\kappa(F)$
denotes the condition number for the sum function $F(\cdot)$. 

In our recent work, we proposed the \emph{Curvature-aided Incremental Aggregated Gradient} (CIAG)~\cite{CIAG} method, which utilizes curvature (\ie Hessian) information to acquire accelerated convergence. Under the asynchronous model of computation~\cite{bertsekas1989parallel} (where each component function is selected by the algorithm infinitely often, with a uniformly bounded delay), CIAG was shown to converge linearly at a rate that is strictly faster than SAGA;
see also~\cite{gower2017tracking} for a similar technique to variance
reduction in stochastic gradient methods.  
Nevertheless, CIAG is not directly amenable to a distributed implementation, 
where it would require message-passing among agents along a closed spanning walk following 
a \emph{fixed} sequence. This is typically 
restrictive, 
especially in distributed wireless sensor networks where time synchronization is challenging~\cite{freris2011fundamental} and it is desirable that agents 
are activated at random time instants 
for the sake of battery lifetimes. 

This paper 
adopts a stochastic setting where the component
functions are selected at random, \ie it allows for a Markov-driven distributed
implementation on a general topology.
Moreover, the proposed \emph{Stochastic Unbiased Curvature-aided Gradient}
(SUCAG) method introduces an \emph{unbiased} 
curvature-aided gradient estimator.  
We analyze the convergence of SUCAG for strongly convex
problems and show that it converges linearly \emph{with high probability} (w.h.p.), for a sufficiently small stepsize selection. 
Besides, 
for an initialization sufficiently close to optimality, the algorithm requires ${\cal O}( \kappa(F) \log (1 / \epsilon))$ 
steps to reach an $\epsilon$-optimal solution,  
\emph{robust to the graph structure or network size}.
This can be beneficial in the prelude of large-scale cyberphysical systems \cite{CPS}
that feature millions of agents, for example in real-time learning from streaming data~\cite{RCS,RCS2}. 
\vspace{.1cm} 

\textbf{Notation}.
We denote vectors and matrices using boldface lower-case letters and upper-case letters, 
respectively. 
We use the notation $(x)_+:=\max\{0,x\}$ for the positive part.
We use the standard Bachmann-Landau notation: for non-negative functions $f,g,h$ we write $f(t) = {\cal O}( g(t) )$ (\resp $f(t) = \Omega(h(t))$) when there exists a positive constant $c$ (\resp $C$) such that $f(t) \leq c g(t)$ ($f(t) \geq C h(t)$). 
We use $\| \cdot \|$ for the standard Euclidean norm.

\subsection{Overview of Gradient Tracking Techniques} \label{sec:ciag}
Before we introduce the proposed SUCAG algorithm, we briefly 
review the curvature-aided gradient tracking technique. In particular, we cast it as a specific variance reduction method 
within the stochastic gradient paradigm. 

When the number of components $N$ is large, 
and especially in a distributed setting, 
the total gradient $\grd F$ is costly to evaluate,  
whence a stochastic gradient (SG) scheme is  
usually sought to tackle \eqref{eq:opt} in a tractable manner. 
To this end, a SG scheme 
can be described as follows: at iteration $k \in \NN$, 
a \emph{stochastic gradient surrogate} ${\bm g}^k$ is obtained, satisfying
$\EE [ {\bm g}^k ] = \grd F( \prm^k ),~{\rm Var}( {\bm g}^k ) \leq \sigma^2$, 
\ie an unbiased estimator for $\grd F(\prm^k)$ with 
finite error variance. A common implementation 
selects one of the component functions uniformly at random, \ie 
the $i_k-$th function (where $i_k \in [N]$) 
at iteration $k$, and evaluates the gradient $\grd f_{i_k}(\prm^k).$ 
 
The SG method~\cite{robbins1951stochastic} then solves \eqref{eq:opt} using the recursion 
$\prm^{k+1} 
= \prm^k - \gamma_k {\bm g}_{\sf SG}^k$  
with ${\bm g}^k_{{\sf SG}}\eqdef\grd f_{i_k}(\prm^k)$, where 
$\gamma_k >0$ is the stepsize at the $k-$th iteration:  when $\{\gamma_k\}$ satisfies 
$\sum_k \gamma_k = \infty, \sum_k \gamma_k^2 < \infty$,  
the SG method converges to an optimal solution of \eqref{eq:opt} 
almost surely (a.s.) \cite{robbins1951stochastic}. 
However, albeit its simplicity, 
a major drawback of SG remains that 
the (expected) convergence rate is \emph{sublinear}, \ie
$\EE [\| \prm^\star - \prm^k \| ]= {\cal O}(1/k)$ (where $\prm^\star$
denotes the optimal solution to \eqref{eq:opt})
when the objective function is strongly convex \cite{moulines2011non}. 

On the other hand, one may adopt an \emph{incremental} approach to the problem. 
In this purview, a reasonable remedy to the high cost of evaluating the total gradient is to \emph{aggregate} previous gradient component evaluations, \ie approximate 
$\grd f_j (\prm^k ) \approx \grd f_j ( \prm^{\tau_{j}^{k}} )$ 
where $\tau_{j}^{k}$ denotes the iteration count when the $j-$th gradient was last evaluated (see \eqref{eq:tau} for 
the precise definition of $\tau_i^k$). This is precisely the mechanism of SAGA~\cite{SAGA} that selects:\vspace{-.2cm}
\beq\label{eq:saga} \notag
{\bm g}^k_{\sf SAGA}\eqdef \left[\grd f_{i_k}(\prm^k) - \grd f_{i_k}(\prm^{\tau_{i_k}^{k-1}})\right] + \frac{1}{N}\sum_{j=1}^N\grd f_{j}(\prm^{\tau_j^{k-1}}).
\vspace{-.1cm}
\eeq 
Besides, SAG is different only in that it multiplies the first term above by $\frac{1}{N}$ and is, therefore, biased.  
When each component function $f_j$ is smooth, the approximation error, \ie the \emph{variance} of 
the total gradient estimator, 
is bounded by 
${\cal O} ( \sum_{i=1}^N \| \prm^{\tau_i^{k}} - \prm^\star \|)$.
Thanks to the \emph{variance reduction} property, the SAGA/SAG  methods  
are shown to converge at a linear rate~\cite{SAGA,SAG} and are therefore faster than SG, at the cost of storing the past gradient component values. It is worthwhile observing that in a multi-agent setting, past values can be stored locally, and aggregated in a recursive manner via message passing. The crucial implementation detail in our algorithm is that aggregation can be performed via simple \emph{gossiping}, \ie message exchange between an agent and only one of its neighbors.

Under an additional smoothness condition (of Lipschitz-continuous component Hessians),  
our recent work \cite{CIAG} considered 
a first-order Taylor approximation to the
total gradient 
 in order to further reduce the estimation variance: \vspace{-.1cm}
\beq \label{eq:ciag}
{\bm g}^k_{\sf CIAG}\eqdef \frac{1}{N}\sum_{j=1}^N \left[\grd f_j( \prm^{\tau_j^{k}} ) + \grd^2 f_j( \prm^{\tau_j^{k}} )
( \prm^k -  \prm^{\tau_j^{k}} )\right],\notag\vspace{-.1cm}
\eeq
A distinctive attribute is that the approximation error for ${\bm g}^k_{\sf CIAG}$ 
to the total gradient 
above can be bounded
by ${\cal O}( \sum_{j=1}^N\| \prm^k - \prm^{\tau_j^k} \|^2 )$ 
in light of~\cite[Lemma 1.2.4]{nesterov_lectures}, 
\ie the bound features \emph{squared norm} 
of the 
distance to optimality; specifically, when $\| \prm^k - \prm^{\tau_j^k} \|^2 < 1$ (as is the case for large enough $k$), this bound allows a tighter convergence analysis as compared to SAG/SAGA. 
In brief, it was shown in \cite{CIAG} 
that CIAG method exhibits a faster linear convergence rate than SAG \emph{by a factor of $N$}, under the restriction that the delays in gradient/Hessian evaluation are
bounded (\ie $| \tau_i^k - k | \leq K < \infty$, e.g., when $i_k$ 
is chosen as $i_k = ( k ~{\rm mod}~ N ) + 1$). Note, however, that in a stochastic setup CIAG is biased, which motivates the development of an unbiased variant in this paper.

\section{The SUCAG method}
We first introduce SUCAG method as a stochastic gradient algorithm for finite sum optimization \eqref{eq:opt}, and  
then proceed to present 
protocols for implementing it in a distributed setting. 

We adhere to the stochastic optimization setting: 
at iteration $k$, the algorithm selects 
the $i_k-$th component function $f_{i_k} (\cdot)$, where $i_k \in \{1,...,N\}$, and evaluates its gradient/Hessian at the current iterate. At first, we assume that the selection mechanism is independent identically distributed (i.i.d.), where a  
random integer is chosen \emph{uniformly at random}, \ie $P(i_k = j) = 1/N$ for all $j =1,...,N$.
We define the random variable $\tau_j^k$ as the iteration number where the $j-$th function was last accessed
(including iteration $k$ itself):
\beq \label{eq:tau}
\tau_j^k := \max\{\ell \geq 0 ~:~ i_\ell = j,~\ell \leq k \} \eqs.
\eeq
Note that $\tau_{i_k}^k = k$ and $\tau_j^k = \tau_j^{k-1}$ for all $j \neq i_k$.

Inspired by the curvature-aided gradient tracking technique~\cite{CIAG},  
the SUCAG method adopts the update rule:
\beq \label{eq:sucag}
\prm^{k+1} = \prm^k - \gamma {\bm g}_{\sf SUCAG}^k,~k \geq 0 \eqs,
\eeq
where $\gamma > 0$ is a fixed step size
and ${\bm g}_{\sf SUCAG}^k$ is a curvature-aided approximation 
on the actual gradient
$\grd F( \prm^k )$. The approximation is given as\footnote{Note that CIAG pre-multiplies the first term with $\frac{1}{N}$, in the same way that SAG does for SAGA.}:
\beq \label{eq:ucag_g}
\begin{split}
& \textstyle {\bm g}_{\sf SUCAG}^k \eqdef \left[ \grd f_{i_k}(\prm^k) - {\bm g}_{i_k}^k (\prm^k) \right] + \frac{1}{N} \sum_{i=1}^N {\bm g}_i^k (\prm^k)  \eqs,
\end{split}
\eeq
where 
\beq \label{eq:sucag_e}
{\bm g}_i^{k} (\prm) \eqdef \grd f_i (\prm^{\tau_i^{k-1}}) + \grd^2 f_i (\prm^{\tau_i^{k-1}}) \big( \prm - \prm^{\tau_i^{k-1}} \big) \eqs.
\eeq
A key property of ${\bm g}_{\sf SUCAG}^k$ is that it is an 
\emph{unbiased} estimator for $\grd F(\prm^k)$. 
To see this, define the filtration ${\cal F}_k$ as the collection of random variables 
(the chosen indices $\{i_k\}$) 
realized prior to the update of $\prm^{k+1}$, \ie ${\cal F}_k\eqdef \sigma(\{i_{\ell}:\ell = 0,1,\hdots,k-1\})$;  
note that this does not include the newly selected index  $i_k$.
The conditional expectation is given by: 
\beq \label{eq:unbiased} \begin{split} 
\EE_k[ {\bm g}_{\sf SUCAG}^k ] &   = \textstyle
\frac{1}{N} \sum_{i=1}^N \big( {\bm g}_i^k (\prm^k) + \grd f_i(\prm^k) - 
{\bm g}_i^k( \prm^k ) \big)  \\
& = \grd F( \prm^k ) \eqs,
\end{split}
\eeq
where $\EE_k[\cdot] \eqdef \EE[ \cdot | {\cal F}_k ]$, and
therefore SUCAG is unbiased. 
Note that SUCAG is a natural extension of SAGA to incorporate curvature information.  

\algsetup{indent=0.8em}
\begin{algorithm}[t]
\caption{SUCAG method}\label{alg:sucag}
  \begin{algorithmic}[1]
  \STATE \textbf{Input}: Initial point $\prm^0 \in \RR^d$
  \STATE Set counter variables $\tau_i^0 \eqdef -1$ for all $i$.
  \STATE Set $\grd f_i( \prm^{-1} ) = {\bm 0}$ and $\grd^2 f_i ( \prm^{-1} ) = {\bm 0}$ for all $i$.
  \FOR {$k=0,1,2,\hdots, K$}
  \STATE Choose $i_k \in \{1,\hdots,N\}$ uniformly at random. 
  \STATE Compute ${\bm g}_{\sf SUCAG}^k$ using \eqref{eq:ucag_g}, 
  \eqref{eq:sucag_e}. 
   \STATE 
   \textbf{Update:}
   \[
   \prm^{k+1} = \prm^k - \gamma {\bm g}^k_{\sf SUCAG} \eqs. \vspace{-.5cm}
   \]
   \STATE Update counter variables: $\tau_{i_k}^k \leftarrow k$, 
   and $\tau_j^k \leftarrow \tau_j^{k-1}$ for all $j \neq i_k$. 
\ENDFOR
\STATE \textbf{Return} 
$\prm^{K+1}$.
  \end{algorithmic}
\end{algorithm}
While our main focus is on applying SUCAG in a distributed setting, we note that one can also implement the 
SUCAG method as a centralized incremental method just as in CIAG~\cite{CIAG}. In this case, the system stores ${\cal O}(Nd)$ real numbers in the memory for the 
previous iterates and the per-iteration
complexity is ${\cal O}(d^2)$, cf.~\cite{CIAG}. 

\subsection{Distributed Implementation}
\label{sec:dist}
Consider a \emph{connected} undirected 
graph $G = (V,E)$ where $V = \{1,\hdots,N\}$ 
is the set of $N$ agents and $E \subseteq V \times V$ is the edge set that prescribes the 
agent-to-agent communication pairs. 
The neighborhood set of $i$ is $\mathcal{N}_i\eqdef \{j: (i,j)\in E\}$. 

It is easy to verify that:
\beq \label{eq:sucag_a}
{\bm g}_{\sf SUCAG}^k = \left[ \grd f_{i_k}(\prm^k) - {\bm g}_{i_k}^k (\prm^k) \right]+ {\bm b}^{k-1} + {\bm H}^{k-1} \prm^k, 
\eeq
where
\beq \notag
\begin{split}
& \textstyle {\bm b}^{k-1} \eqdef \frac{1}{N} \sum_{i=1}^N \big( \grd f_i( \prm^{\tau_i^{k-1}} ) - 
\grd^2 f_i( \prm^{\tau_i^{k-1}} ) \prm^{\tau_i^{k-1}} \big) \eqs, \\
& \textstyle {\bm H}^{k-1} \eqdef \frac{1}{N} \sum_{i=1}^N  \grd^2 f_i( \prm^{\tau_i^{k-1}} ) \eqs.
\end{split}
\eeq
These variables can be computed recursively as follows:
\beq \label{eq:inc_upd} 
\begin{split}
& \textstyle {\bm b}^{k} \eqdef {\bm b}^{k-1} + \frac{1}{N} \left( \grd f_{i_k}( \prm^k ) - \grd f_{i_k}( \prm^{\tau_{i_k}^{k-1}})\right) \\
& \textstyle \hspace{1.35cm} \frac{1}{N}\left(\grd^2 f_i(  \prm^{\tau_{i_k}^{k-1}} )  \prm^{\tau_{i_k}^{k-1}}- \grd^2 f_i( \prm^k ) \prm^k \right)\\
& \textstyle {\bm H}^{k} \eqdef {\bm H}^{k-1} + \frac{1}{N} \big(  
\grd^2 f_{i_k} ( \prm^k ) -  \grd^2 f_{i_k} ( \prm^{\tau_{i_k}^{k-1}} ) 
\big) \eqs.
\end{split}
\eeq
Leveraging these recursions, we can obtain distributed protocols for implementing SUCAG on $G$,  
via transmitting the aggregated (previous) gradients and Hessians. To this end, a key aspect is that the required updates for SUCAG [cf.~\eqref{eq:sucag}, \eqref{eq:sucag_a}, \eqref{eq:inc_upd}] are \emph{locally computable} by an activated agent $i_k$, given the latest estimate $\prm^k$ and aggregate gradient/Hessian information ${\bm b}^{k-1}, {\bm H}^{k-1}$. 

First, for the special case when $G$ is a \emph{star graph}, a simple protocol can be implemented by using the hub node as a coordinator to store the tuple $\{{\bm b}^{k-1}, {\bm H}^{k-1}, \prm^k\}$ and, at each iteration: \emph{(i)} the coordinator activates an agent independently and uniformly at random and passes along this tuple, \emph{(ii)} the agent performs the required updates [cf.~\eqref{eq:inc_upd},\eqref{eq:sucag_a},\eqref{eq:sucag}] and \emph{(iii)} transmits back the updated values $\{{\bm b}^k, {\bm H}^k, \prm^{k+1}\}$. 
Such protocol is ideal for cloud computing services.

\begin{figure}[t]
\centering
\ifplainver
{\sf \resizebox{.6\linewidth}{!}{\input{./Tikz/SUCAG.tikz.tex}}}\vspace{-.1cm}
\else
{\sf \resizebox{.9\linewidth}{!}{\input{./Tikz/SUCAG.tikz.tex}}}\vspace{-.1cm}
\fi
\caption{Distributed implementation of SUCAG; the agent sequence $\{ i_k \}_{k \geq 1}$ corresponds to a random walk on the graph $G$.}\vspace{-.4cm} \label{fig:sucag}
\end{figure}
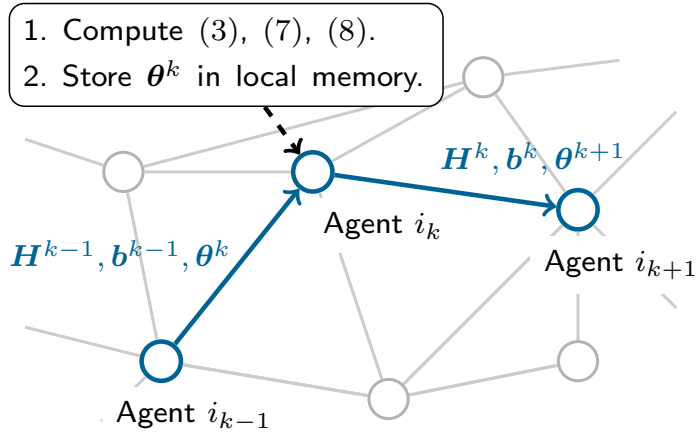

A crucial point of our analysis (cf.~Theorem~\ref{thm:main}) 
reveals that the requirement of independent and uniform agent activation \emph{can be relaxed}. 
In fact, it is adopted primarily for demonstration purposes 
(to ensure unbiasedness, $\EE_k[ {\bm g}^k_{\sf SUCAG} ] = \grd F(\prm^k)$) 
so as to portray our scheme within the popular stochastic variance 
reduction framework. 
This important attribute allows for distributed computations along a random walk on
a connected graph $G$ as follows: 
at iteration $k$, activated agent $i_k$ performs  the update~\eqref{eq:sucag} and \emph{(i)} stores $\prm^{k}$\footnote{In light of storage/computation trade-off, it is also possible for agent $i_k$ to store $\{\prm^{k+1},\grd f_{i_k}(\prm^k),\grd^2f_{i_k}(\prm^k)\}$ in order to avoid re-evaluation at the next activation instance.}, 
\emph{(ii)} chooses one of its neighbors $i_{k+1} \in {\cal N}_{i_k}$ at random 
as the next active agent, and \emph{(iii)} transmits the updated $\{{\bm b}^k, {\bm H}^k, \prm^{k+1}\}$ 
to agent $i_{k+1}$. This process is illustrated in 
Fig.~\ref{fig:sucag}. 

We remark that the above line of reasoning directly applies to \emph{any} stochastic algorithm with the \emph{random selection-aggregate-update} paradigm, e.g., SG, SAG and SAGA; see \cite{johansson2009randomized,ram2009incremental}.

\section{Convergence Analysis}
This section establishes 
the linear convergence of the SUCAG method 
with high probability. 
We begin by stating the required assumptions
on the objective function of \eqref{eq:opt}:\vspace{-.2cm}
\begin{Assumption}\label{ass:cts_h}
Each component function $f_i (\prm)$ has $L_{H,i}$-Lipschitz continuous Hessian. In other words, for any $i\in \{1,...,N\}$, there exists $L_{H,i}>0$, such that for all $\prm', \prm \in \RR^d$ \vspace{-.1cm}
\beq \begin{split}
& \| \grd^2 f_i ( \prm ) - \grd^2 f_i (\prm' ) \| \leq L_{H,i} \| \prm - \prm' \| \eqs.
\end{split} \vspace{-.1cm}
\eeq
We define $\bar{L}_H \eqdef \max_{i\in \{1,...,N\}} L_{H,i}$.
\end{Assumption}\vspace{-.4cm}
\begin{Assumption}\label{ass:cts}
The gradient of sum function $F(\prm)$ is $L$-Lipschitz continuous, 
\ie 
there exists $L>0$, such that for all $\prm', \prm \in \RR^d$, \vspace{-.4cm}
\beq \begin{split}
& \| \grd F ( \prm ) - \grd F (\prm' ) \| \leq L \| \prm - \prm' \| \eqs. 
\end{split}
\eeq
\end{Assumption}\vspace{-.4cm}
\begin{Assumption}\label{ass:strcvx}
The sum function $F (\prm)$ is $\mu$-strongly convex, $\mu > 0$, \ie for all $\prm', \prm \in \RR^d$, 
\beq
F( \prm' ) \geq F( \prm ) + \langle \grd F(\prm), \prm' - \prm \rangle + \frac{\mu}{2} \| \prm' - \prm \|^2 \eqs.
\eeq
\end{Assumption}
Under A\ref{ass:strcvx}, a unique optimal solution 
to problem~\eqref{eq:opt} exists and it is denoted by $\prm^\star$. 
We denote the condition number for $F(\prm)$ 
by $\kappa(F) \eqdef L / \mu$. 
Note that the above are standard assumptions that can be satisfied in many 
applications; for example, $\ell_2-$regularized logistic regression~\cite{vapnik1999overview}, and quadratic programming such as (overdetermined) least-squares and Model Predictive Control (MPC)~\cite{morari1999model} applications.

We also state the following assumption on the delays in the SUCAG method:\vspace{-.2cm}
\begin{Assumption} \label{ass:delay} \label{ass:delay2}
	For all $k \geq 1$, it holds that $|k - \tau_i^{k-1}| \leq m_k$ for all $i$, where $\{m_k\}_{k \geq 1}$ is non-decreasing with $m_1\ge 1$. Moreover, it holds that \vspace{-.1cm}
\beq \label{eq:delay}
4 \log m_k  \leq c_0 + \log  k
\eqs, \vspace{-.1cm}
\eeq
for all $k \geq 1$ and some $c_0\ge0$.
\end{Assumption}
Note that it follows that A\ref{ass:delay2} implies that $\{m_k\}_{k \geq 1}$ 
satisfies \vspace{-.1cm}
\beq 
2 m_k \leq m_0 + 
\big( (1/3) - \beta \big) k, \vspace{-.1cm}
\eeq
for proper selection of $m_0>0$ and $\beta\in (0,\frac{1}{3})$. 
In the interest of space, the detailed proofs
for the below analysis can be found in our 
online appendix (\texttt{https://arxiv.org/abs/1803.08198}).

Later, we will show that the distributed message passing 
protocols in Section~\ref{sec:dist} satisfies A\ref{ass:delay} with high probability. 
The following establishes
the linear convergence of SUCAG:
\begin{Theorem} \label{thm:main}
Assume that A\ref{ass:cts_h}, A\ref{ass:cts}, A\ref{ass:strcvx}, A\ref{ass:delay} 
hold, and let\footnote{This is with absolutely no loss in generality, and is solely used in the analysis so as to simply notation in the derived stepsize selection rule.} $\mu \geq 1$. 
Fix $\delta = 1 - \Delta \gamma \frac{{2}\mu L}{\mu+L}$ for some sufficiently small $\Delta \in (0,1)$, 
and step size $\gamma$ that satisfies:
\beq \label{eq:gamma}
\begin{split}
& \gamma \leq \frac{ 2 }{ \mu + L}\\
& \gamma \leq \min_{j=1,...,4} \Big( \frac{3m_0 }{2} \frac{\mu L}{\mu + L} \Delta + \frac{1}{ C_j'} \Big)^{-1} \eqs,
\end{split}
\eeq
where the constants are given by:
\beq
C_j' \eqdef \frac{e^{1-c_0} \beta \Delta(1-\Delta) }{4 \tilde{q}_j \| \prm^0 - \prm^\star \|^{2 \tilde\eta_j-2} } \Big( \frac{2 \mu L}{\mu + L} \Big)^2,~j =1,...,4 \eqs,
\eeq
with
\beq \notag
\begin{split}
& \tilde{q}_1 \eqdef 2 \bar{L}_H L^2, \tilde{q}_2 \eqdef 32 \bar{L}_H^3, \tilde{q}_3 \eqdef  8 \bar{L}_H^2 L^4, \tilde{q}_4 \eqdef 2048 \bar{L}_H^6, \\
& \tilde\eta_1 \eqdef 3/2,~\tilde\eta_2 \eqdef 5/2,~\tilde\eta_3 \eqdef 2,~\tilde\eta_4 \eqdef 4 \eqs.
\end{split}
\eeq
Then, the  following properties hold:
\begin{enumerate}[leftmargin=.5cm]
\item The SUCAG 
method converges linearly as $ \| \prm^k - \prm^\star \|^2 \le \delta^k \| \prm^0 - \prm^\star \|^2$
for all $k \geq 0$.
\item There exists an upper bound sequence $\{ \bar{V}(k) \}_{k \geq 0}$ 
such that $\| \prm^k - \prm^\star \|^2 \leq \bar{V}(k)$ for all $k \geq 0$ with
the rate:
\beq \label{eq:asympt}
\lim_{ k \rightarrow \infty} \frac{ \bar{V}(k+1) } {\bar{V}(k) } = 1 - 2 \gamma \frac{ \mu L }{ \mu +L} \eqs.
\eeq
\end{enumerate}
\end{Theorem}
First, 
from Theorem~\ref{thm:main}, 
we consider the case when 
$\| \prm^0 - \prm^\star \| \approx 0$.
In particular, taking 
$\Delta = \| \prm^0 - \prm^\star\|^{\frac{1}{2}}$ shows that 
the right hand side
of the second inequality in \eqref{eq:gamma} 
is at the order of $\Omega( \| \prm^0 - \prm^\star \|^{-\frac{1}{2}} )$. 
When $\| \prm^0 - \prm^\star \| \approx 0$, the SUCAG 
method is allowed to take the step size at $\gamma = 2 / (\mu+L)$
in \eqref{eq:gamma}. 
Substituting this into \eqref{eq:asympt} shows that
the asymptotic linear convergence rate of the SUCAG method
is $1 - 4 \mu L / (\mu + L)^2 = 1 -\Omega( 1/\kappa(F))$. 
Consequently, the method finds an $\epsilon$-optimal solution
to \eqref{eq:opt} using just ${\cal O}( \kappa(F) \log(1/\epsilon) )$
steps. 

Second, assumption A\ref{ass:delay} on the delays in the SUCAG method 
allows for a sub-linear delay which 
grows \emph{unboundedly} as fast as $m_k = {\cal O}(k^{1/4})$. 
This is of relevance to the two 
distributed implementations discussed in Section~\ref{sec:dist}. 
In particular, in both cases
the agents' activations  
can be captured by a finite state Markov chain, and the time delays can be
analyzed via the \emph{hitting time} of the chain:\vspace{-.1cm}
\begin{Lemma} \label{lem:mc}
Consider the sequence $( i_k )_{k \geq 0}$ 
that is governed by a finite, irreducible and aperiodic Markov chain, and suppose that the chain is 
stationary\footnote{We note that the analysis applies even to non stationary Markov chains; we skip the details for length considerations.}. 
It holds that
\hide{
\beq \label{eq:mc}
P( | k - \tau_i^{k-1} | > 1+x ) \leq {\rm exp} \Big( -\frac{x}{1 + e \bar{\tau}_i} + 1 \Big),~\forall~x \geq 0,
\eeq 
}
\beq \label{eq:mc}
P( | k - \tau_i^{k-1} | > x ) \leq {\rm exp} \Big( -\frac{x}{1 + e \bar{\tau}_i} + 1 \Big),~\forall~x \geq 0,
\eeq 
for all $i =1,...,N$, where $\bar{\tau}_i$ is the expected first visit time. 
In particular, the  delay
$| k - \tau_i^{k-1}|$ can be bounded by ${\cal O}( \bar{\tau}_i + \bar{\tau}_i \log (Nk^2/\epsilon)$ with probability at least $1-\epsilon$.
\end{Lemma}
In specific, for any connected undirected graph, 
the Markov chain is stationary and $\bar{\tau}_i = N$ 
provided that the first activated agent is selected uniformly at random. 
Consequently, even though $m_k \rightarrow \infty$ as $k \rightarrow \infty$, the time varying delay $m_k$ satisfies assumption A\ref{ass:delay}
w.h.p.

Importantly, combining Theorem~\ref{thm:main} and Lemma~\ref{lem:mc} 
shows that when the initialization is sufficiently
close to the optimal solution, w.h.p.~the distributed implementation of SUCAG converges linearly at an asymptotic rate that is \emph{independent from the communication graph structure}.

\subsection{Proof of Theorem~\ref{thm:main}} \label{sec:pf}
The proof consists of two parts --- first we prove a descent lemma for the SUCAG 
method to characterize its per-iteration progress; then we show that the
resulting nonlinear process converges linearly, with asymptotic rate as in \eqref{eq:asympt}. 
Let us state the first result as follows:\vspace{-.1cm}
\begin{Lemma} \label{lem:descent}
Assume that A\ref{ass:cts_h}, A\ref{ass:cts}, A\ref{ass:strcvx}, A\ref{ass:delay} hold
and $\gamma \leq 2 / (\mu+L)$. 
The following holds for all $k \geq 0$:
\beq \label{eq:lemma1}
\begin{split}
& \| \prm^{k+1} - \prm^\star \|^2 \leq \Big( 1 - 2 \gamma \frac{\mu L}{\mu + L}\Big) \| \prm^k - \prm^\star \|^2 \\
& +  \gamma^3 (m_k)^2  \max_{ (k-2m_k) \leq \ell \leq k } \Big( \tilde{q}_1 \| \prm^{\ell} - \prm^\star \|^3 + 
\tilde{q}_2 \| \prm^{\ell} - \prm^\star \|^5 \Big) \\
 & + \gamma^6 (m_k)^4 \max_{ (k-2m_k) \leq \ell \leq k } \Big( 
 \tilde{q}_3 \| \prm^\ell - \prm^\star \|^4 + 
 \tilde{q}_4 \| \prm^\ell - \prm^\star \|^8 \Big)
 \\[-.2cm] \end{split}
\eeq 
\end{Lemma}
An important fact observed is that the 
distance to optimal solution at the $(k+1)-$th iteration 
follows an inequality system  with high order terms 
of the distance to optimal solution evaluated at \emph{delayed} time instances. 

Let $R(k) \eqdef \| \prm^k - \prm^\star \|^2$, 
the above motivates us to study the system:
\beq \label{eq:lem_ineq}
R(k+1) \leq p R(k) + \sum_{j=1}^J q_j  m_{1,j}^{(k)} \hspace{-.1cm} \max_{ (k - m_{2,j}^{(k)})_+ \leq \ell \leq k } \hspace{-.1cm}
R^{\eta_j}(\ell) , 
\eeq
where $p\in (0,1)$, $\eta_j > 1$ for all $j$ and $R(k)$ is non-negative. 
Moreover, we have $m_{1,j}^{(k)}, m_{2,j}^{(k)} \geq 0$.
The following lemma establishes the linear convergence for $R(k)$
under a sufficient condition: 
\begin{Lemma} \label{lem:nonlinear}
Fix $\delta \in (0,1)$  such that $\delta > p$. 
Assume that 
\beq \label{eq:delay_gen}
\log( m_{1,j}^{(k)} ) \leq c_{0} + c_{1} \log k,~
m_{2,j}^{(k)} \leq m_0 + \Big( \frac{\eta-1}{\eta} - \beta \Big) k \eqs,
\eeq
for some $m_0, c_{0}, c_{1} \geq 1$, $0 < \beta < (\eta-1)/\eta$, 
and $\eta \eqdef \min_{j=1,...,J} \eta_j > 1$.
Moreover, assume that
\beq \label{eq:lem_1_pre}
q_j \leq \frac{1}{R^{\eta_j-1}(0)} \frac{ \delta^{-{\xi}_j^\star(\delta)} (\delta - p) }{J} \eqs,~j=1,\hdots,J \eqs,
\eeq
with the non-positive number $\xi_j^\star(\delta)$ defined as:
\beq \label{eq:lem_2}
\xi_j^\star ( \delta) \eqdef \min_{ k \geq 0} \Big( \frac{ \log m_{1,j}^{(k)} }{\log \delta} + \eta_j \big( k - m_{2,j}^{(k)} \big)_+ - k \Big) \eqs.
\eeq
It holds for \eqref{eq:lem_ineq} that \emph{(i)} $R(k) \leq \delta^k R(0)$ for all $k \geq 0$, 
and \emph{(ii)} there exists an upper bound sequence $\{ \bar{R}(k) \}_{k \geq 0}$ such that $\bar{R}(k) \geq R(k)$ for all $k \geq 0$ and its convergence rate is:
\beq
\lim_{ k \rightarrow \infty } \bar{R}(k+1) / \bar{R}(k) = p \eqs.
\eeq
\end{Lemma}
We remark that a crucial point to guaranteeing the linear convergence rate 
is on the high order terms with power $\eta_j > 1$. 
In fact, for $\eta_j = 1$ with some $j$, it can be shown 
that $R(k)$ may diverge as $m_{1,j}^{(k)}, m_{2,j}^{(k)}$ is growing sublinearly in \eqref{eq:delay_gen}. 

Finally, 
we observe that \eqref{eq:lemma1} is a 
special case of \eqref{eq:lem_ineq}
by substituting $R(k) = \| \prm^k - \prm^\star \|^2$ and the parameters:
\beq \nonumber
\begin{split}
& p = 1 - \frac{2 \gamma \mu L}{\mu+L},~\eta_1 = \frac{3}{2},~\eta_2 = \frac{5}{2},~\eta_3 = 2,~\eta_4 = 4 \eqs, \\
& q_1 = \gamma^3 \tilde{q}_1,~q_2 = \gamma^3 \tilde{q}_2,~q_3 = \gamma^6 \tilde{q}_3,~q_4 = \gamma^6 \tilde{q}_4 \eqs,\\
& m_{1,j}^{(k)} = (m_k)^{2 \lceil \frac{j}{2} \rceil},~m_{2,j}^{(k)} = 2m_k,~j=1,\hdots,4  \eqs,
\end{split}
\eeq
Under A\ref{ass:delay}, we observe that the above $m_{1,j}^{(k)}, m_{2,j}^{(k)}$ 
satisfy \eqref{eq:delay_gen} in 
Lemma~\ref{lem:nonlinear} with the same $c_0$ . 
Moreover, in the online appendix we show 
that the step size $\gamma$ specified in \eqref{eq:gamma} 
gives the set of $q_j$ that satisfy \eqref{eq:lem_1_pre}. 
This concludes the proof.

\hide{
Bound that we use:
\nfc{c/p the analysis for uniform distribution}

\nfc{we also need a small lemma for MC large deviations: the two papers we found so far are in Dropbox: ToRead$\rightarrow$MC\_Chernoff}.}

\section{Numerical Experiments}
\begin{figure}[t]
\centering
\ifplainver
\resizebox{.6\linewidth}{!}{\input{./Tikz/cdc_test_case.tikz.tex}}
\else
\resizebox{.9\linewidth}{!}{\input{./Tikz/cdc_test_case.tikz.tex}}
\fi
\caption{Comparison of distributed optimization methods for solving the logistic regression problem \eqref{eq:logist} over a time varying graph
with $N=250$ agents, batch size $B=1$, and $d=51$ parameters. The optimality gap 
against iteration number is obtained through averaging over $10$ trials (except for the deterministic CIAG method).}\vspace{-.3cm} \label{fig:sim}
\end{figure}
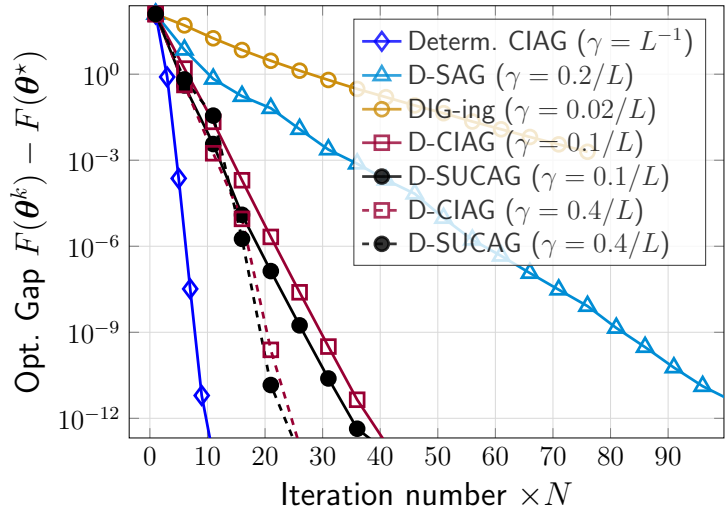

We evaluate the performance SUCAG on the learning task of training a linear classifier via logistic regression. In this problem, the $i-$th 
component function is
defined as:
\beq \label{eq:logist}
f_i ( \prm ) \eqdef \sum_{b=1}^B \log( 1 + {\rm exp} ( -y_{i,b} \langle \prm, {\bm x}_{i,b} \rangle ) ) + \frac{B}{2N}\| \prm \|^2 \eqs,
\eeq
where $( {\bm x}_{i,b}, {\bm y}_{i,b} )_{b=1}^B$ represents  
data held by agent $i$, and $B$ is the size of the corresponding data subset.  
It can be verified that \eqref{eq:logist} and its corresponding sum function $F(\prm)$ 
satisfy A\ref{ass:cts_h}--A\ref{ass:strcvx}. 
We generate a synthetic dataset in the experiment 
where the ground truth classifier $\prm_0$ is selected as
$\prm_0 \sim {\cal U}[-1,1]^d$, then each feature vector ${\bm x}_i$ 
is generated independently following ${\bm x}_i \sim {\cal U}[-1,1]^{d}$ 
with the label $y_i$ computed as $y_i = {\rm sign} ( \langle {\bm x}_i, \prm_0 \rangle )$.
To simulate the distributed optimization methods, we generate $G$ as an Erdos-Renyi graph with $N$ agents and a connection probability of $2 \frac{\log N}{N}$ (where, if needed, we repeat the procedure until a connected graph is obtained). 
The SUCAG method operates on the aforementioned graph using the protocol
in Section~\ref{sec:dist} (see also~Fig.\ref{fig:sucag}), where the probability 
for selecting the next agent is uniform over all neighbors.  
For benchmark purposes, we compare the SUCAG method to the similarly
modified version of distributed CIAG~\cite{CIAG} and SAG~\cite{SAG}.  
Moreover, we compare the DIG-ing method \cite{nedic2017achieving},  
an average consensus-based distributed optimization method that accounts for time varying communication topologies. For this method, the mixing matrix is designed in a similar fashion as the  pairwise gossip protocol \cite{boyd2006randomized}, while the  agents' activation sequence follows the same model  as in our proposed scheme.
The deterministic CIAG method \cite{CIAG} implemented on a fixed ring graph
is also compared for benchmarking purpose.

Fig.~\ref{fig:sim} shows 
the numerical findings (averaged over 10 generated random walks) for a synthetic dataset example 
with $d=51$ parameters, $N=250$ agents, where each 
agent has a batch of $B = 1$ data tuple, \ie
there are $NB = 250$ data tuples in total. 
For both SAG and DIG-ing methods, we have optimized their
step sizes (as $\gamma = 0.2/L$ and $\gamma = 0.02/L$, respectively) 
used so that the methods converge in the trials 
performed. 
For SUCAG and distributed CIAG methods, two step size configurations are compared.
The numerical results show that the SUCAG method (as well as the distributed CIAG
method) outperform all other methods, and thus clearly demonstrating the benefits of 
curvature-aided gradient tracking 
in accelerating distributed optimization. 

\section{Conclusions}
We have proposed an unbiased curvature-aided stochastic gradient method for large-scale optimization.  
Our method is implementable in a distributed setting, both via distributed computations performed from agents activated by a single coordinator, as well as via a fully distributed random walk on the communication graph. 
Under usual assumptions, our convergence analysis establishes linear convergence at a rate that solely depends on the condition number, but is \emph{independent from the network size} provided that the initialization is sufficiently close to optimality. Besides, we have established more general results under rather weak requirements, \ie   infinitely often activation, even when activation delays grow unbounded: our line of reasoning  may be applied to a wide class of stochastic optimization primitives. Our numerical experiments verified faster convergence compared to various methods for distributed and stochastic optimization.\vspace{.1cm}

\textbf{Acknowledgements}.
This work was supported by NSF under grants NSF CCF-BSF 1714672, CCF-1717207 and CCF-1717391. The authors would like to thank Prof. Maxim Raginsky for pointing out the reference~\cite{roch2015}.  

\bibliography{references}

\ifconfver

\else
\appendix

\ifplainver
\section{Proof of Lemma~\ref{lem:mc}}\label{app:mc}
\else
\subsection{Proof of Lemma~\ref{lem:mc}}\label{app:mc}
\fi
Observe that the random variable $k - \tau_i^{k-1}$ is upper bounded by 
the \emph{return time} for node $i$, for all $k$. 
This is because at time $k$, for $i \neq i_k$, node $i$ has not been visited again since the last time it was visited;  while for $i = i_k$, this quantity is precisely the return time. By ergodicity, it follows that the distribution of the return time equals that for the hitting time whence \cite[Lemma 3.26]{roch2015} applies. 
Upper bounding the right hand side of~\eqref{eq:mc} 
by $\frac{\epsilon}{Nk^2}$, and invoking the union bound concludes the proof. 

\ifplainver
\section{Proof of Lemma~\ref{lem:descent}} \label{app:err}
\else
\subsection{Proof of Lemma~\ref{lem:descent}} \label{app:err}
\fi
Let us define  ${\bm e}^k \eqdef {\bm g}_{\sf SUCAG}^k - \grd F(\prm^k)$.
Under A\ref{ass:cts_h}, A\ref{ass:cts}, A\ref{ass:strcvx}, we 
observe that the simple analysis leading to Eq.~(36), (37) from \cite{CIAG} applies directly\footnote{In particular, this requires that $\gamma < \frac{2}{\mu+L}$, which in light of the requirement that $\gamma \le \frac{\mu +L}{2\mu L\Delta}$ in Theorem~\ref{thm:main}, suggests selecting $\Delta < \frac{(\mu +L)^2}{4\mu L}$.}, 
in particular
\beq \label{eq:first_ineq}
 \begin{split}
 \| \prm^{k+1} - \prm^\star \|^2 & \leq \Big(1 -  2\gamma \frac{ \mu L }{ \mu + L} \Big) 
 \| \prm^k - \prm^\star \|^2 \\
 & + \gamma \| \prm^k - \prm^\star \| \| {\bm e}^k \| + \gamma^2 \| {\bm e}^k \|^2 \eqs.
 \end{split}
\eeq
Our task is to obtain an upper bound on $\|{\bm e}^k\|$. 
Observe that:
\beq \label{eq:e_bound1}
\begin{split}
\| {\bm e}^k \| & \textstyle = \Big\| \frac{1}{N} \sum_{i \neq i_k}
\Big( {\bm g}_i^k ( \prm^k )  - \grd f_i (\prm^k) \Big) \\
& \textstyle \hspace{0.75cm} +\Big(1 - \frac{1}{N} \Big) \Big( \grd f_{i_k} ( \prm^{k} ) - {\bm g}_{i_k}^k ( \prm^k ) \big) \Big) {\Big\|}\\
& \textstyle \leq \frac{1}{2N} \sum_{i \neq i_k} L_{H,i} \Big\| \prm^k - \prm^{\tau_i^k} \|^2 \\
& \textstyle \hspace{0.75cm} + \frac{1}{2} \Big(1-\frac{1}{N} \Big) L_{H,i_k} \| \prm^k - \prm^{\tau_{i_k}^{k-1}} \|^2 \eqs,
\end{split}
\eeq
where $g_i^k( \prm)$ was defined in \eqref{eq:sucag_e},
we have applied triangle inequality and the 
Lipschitz Hessian condition (A\ref{ass:cts_h})~\cite[Lemma 1.2.4]{nesterov_lectures} 
to obtain the second inequality. 
For $i \neq i_k$, 
\beq 
\begin{split}
\| \prm^k - \prm^{\tau_i^k} \|^2 & \textstyle \leq (k-\tau_i^k) \sum_{j=\tau_i^k}^{k-1} \| \prm^{j+1} - \prm^j \|^2 \\
& \leq \textstyle m_k \sum_{j=(k-m_k)_+}^{k-1} \| \prm^{j+1} - \prm^j \|^2 \eqs,
\end{split}
\eeq
where we have used $\tau_i^{k-1} = \tau_i^k$ for $i \neq i_k$
and $|k - \tau_i^{k-1}| \leq m_k$. 
Similarly, we also have 
$\| \prm^k - \prm^{\tau_{i_k}^{k-1}} \|^2 \leq m_k \sum_{j=(k-m_k)_+}^{k-1} \| \prm^{j+1} - \prm^j \|^2$. 
Substituting the above into \eqref{eq:e_bound1}, we get:
\beq \notag
\begin{split}
\| {\bm e}^k \| & \leq \frac{ m_k \bar{L}_H }{2} \Big( \frac{N-1}{N} + 1 - \frac{1}{N} \Big) \hspace{-.2cm} \sum_{j=(k-m_k)_+}^{k-1} \hspace{-.2cm} \| \prm^{j+1} - \prm^j \|^2.
\end{split}
\eeq
Note that $\frac{N-1}{N} + 1 - \frac{1}{N} \leq 2$ and
using the fact that $\prm^{j+1} - \prm^j = - \gamma {\bm g}^k_{\sf SUCAG} = -\gamma ( \grd F (\prm^j) + {\bm e}^j )$
we obtain:
\beq \label{eq:e_bound2}
\| {\bm e}^k \|  \leq 2 m_k  \bar{L}_H \gamma^2 \sum_{j= k - m_k}^{k-1} \Big( \| {\bm e}^j \|^2 + \| \grd F (\prm^j)  \|^2 \Big)
\eeq
which is analogous to (42) in \cite{CIAG}. Alternatively, we can 
upper bound $\| {\bm e}^k \|$ 
as
\beq \notag
\begin{split}
\| {\bm e}^k \| & \leq \frac{\bar{L}_H}{2} \Big( \frac{1}{N} \sum_{i \neq i_k} \| \prm^k - \prm^{\tau_i^k} \|^2 + \frac{N-1}{N} \| \prm^k - \prm^{\tau_{i_k}^{k-1}} \|^2 \Big) \\
& \leq \bar{L}_H \Big( 2 \| \prm^k - \prm^\star \|^2 + \| \prm^{\tau_{i_k}^{k-1}} - \prm^\star \|^2 \\
& \hspace{2cm} \textstyle + \frac{1}{N} \sum_{ i \neq i_k } \| \prm^{\tau_i^{k-1}} - \prm^\star \|^2 \Big) \\
& \leq \bar{L}_H \max_{ \ell \in \{ \{ \tau_i^{k-1} \}_{i =1}^N , k \} } \| \prm^\ell - \prm^\star \|^2
\\
& \leq \bar{L}_H \max_{ ( k - m_k ) \leq \ell \leq k } \| \prm^\ell - \prm^\star \|^2 \eqs.
\end{split}
\eeq
Substituting the above back into \eqref{eq:e_bound2} and using 
$\| \grd F(\prm^j) \|^2 \leq L^2 \| \prm^j - \prm^\star \|^2$ [cf.~A\ref{ass:cts}] gives:
\beq \notag
\begin{split}
\| {\bm e}^k \| & \leq \gamma^2 m_k \hspace{-.3cm}
\sum_{j= k - m_k}^{k-1} \hspace{-.3cm} \Big( q_1 \| \prm^j - \prm^\star \|^2 + q_2 \hspace{-.25cm} \max_{ ( j - m_j ) \leq \ell \leq j } \hspace{-.25cm} \|\prm^\ell - \prm^\star \|^4 \Big) \\
& \hspace{-.5cm} \leq ( \gamma m_k)^2 \hspace{-.2cm} \max_{ (k-2m_k)_+ \leq \ell \leq k } \Big( q_1 \| \prm^\ell - \prm^\star \|^2  + q_2 \| \prm^\ell - \prm^\star\|^4 \Big),
\end{split}
\eeq
where the last inequality is obtained by relaxing the set considered in the maximization
and the fact that $m_{k-m_k} \leq m_k$ as $m_k$ is non-decreasing. 
Obviously, it follows that 
\beq \notag
\| {\bm e}^k \|^2 \leq 
( \gamma m_k)^4 \hspace{-.3cm} \max_{ (k-2m_k)_+ \leq \ell \leq k } \Big( q_3 \| \prm^\ell - \prm^\star \|^4  + q_4 \| \prm^\ell - \prm^\star\|^8 \Big),
\eeq
as we note that $q_3 = 2 (q_1)^2$ and $q_4 = 2 ( q_2)^2$. 
Substituting the above back into \eqref{eq:first_ineq} shows the desired result. 

\ifplainver
\section{Proof of Lemma~\ref{lem:nonlinear}} \label{app:hot}
\else
\subsection{Proof of Lemma~\ref{lem:nonlinear}} \label{app:hot}
\fi
For any $\delta \in (0,1)$, the non-positive number
$\xi_j^\star(\delta)$ is well defined for $m_1^{(k)}, m_2^{(k)}$ satisfying \eqref{eq:delay_gen}.  
To see this, note that by assumption
\beq \label{eq:xi_lb}
\begin{split}
\xi^\star( \delta ) & \geq \min_{ k \geq 0 } \frac{ c_0 + c_1 \log k }{\log \delta} + k \Big( (\eta - 1) - \eta \frac{m_{2,j}^{(k)}}{k} \Big) \\
& \geq \frac{ c_0 + c_1 \log k }{\log \delta} + \beta k - {\eta m_0},
\end{split}
\eeq
where $\log \delta < 0$. The right hand side in the above is a {strictly} convex function 
that attains a {unique} minimum {(over $k\in \mathbb{R}_+$)}.
Consequently, 
condition \eqref{eq:lem_1_pre} implies that
\beq \label{eq:lem_1} 
\begin{split}
& ~\textstyle p + \sum_{j=1}^J q_j R^{\eta-1}(0) \delta^{\underline{\xi}^\star (\delta)} \leq \delta \\
\Longrightarrow &~ \textstyle p + \sum_{j=1}^J q_j R^{\eta-1}(0) \delta^{\xi_j^\star (\delta)} \leq \delta \eqs.
\end{split}
\eeq

From now on, we focus on a simplified case with $J=1$ to highlight on the
mechanics involved, whence the subscripts  
on $\eta_j, q_j, m_{1,j}^{(k)}, m_{2,j}^{(k)}$ are omitted. The general case follows the exact same lines of analysis. 
The first step is to use induction to prove the first statement ---  
consider the base case with $k=1$, 
\beq \notag
R(1) \leq p R(0) + q R(0)^{\eta} m_{1}^{(0)} \leq p R(0) + q R(0)^{\eta} \delta^{\xi_\delta^\star} \leq \delta R(0) \eqs,
\eeq
where the second last inequality follows directly from the definition of $\xi_\delta^\star$.
For the induction step, we assume that $R(k) \leq \delta^k R(0)$ and observe the chain:
\beq \label{eq:lem2}
\begin{split}
& R(k+1) \\
& \leq p R(k) + q m_1^{(k)} \max_{( k- m_2^{(k)} )_+ \leq \ell \leq k } R(\ell)^{\eta} \\
& \leq p \cdot \delta^k R(0) + q m_1^{(k)} \delta^{ \eta (k - m_2^{(k)})_+ } R(0)^\eta \\
& = \delta^k R(0) \cdot \Big( p + q \delta^{ \frac{\log m_1^{(k)}}{\log \delta} + \eta (k-m_2^{(k)})_+ - k} R(0)^{\eta-1}\Big) \\
& \leq \delta^k R(0) \cdot \Big( p + q R(0)^{\eta-1} \delta^{\xi^\star(\delta)} \Big) \leq \delta^{k+1} R(0) \eqs,
\end{split}
\eeq
where we have used \eqref{eq:lem_1} in the last inequality. 
We have thus established $R(k) \leq \delta^k R(0)$ for all $k \geq 1$. 

To prove the second statement, we consider an upper bound sequence 
$\{ \bar{R}(k) \}_{k \geq 0}$ such that $\bar{R}(0) = R(0)$ and 
\beq
\bar{R}(k+1) = p\bar{R}(k) + q m_1^{(k)} \max_{ (k-m_2^{(k)} \leq \ell \leq k } \bar{R}^{\eta}(\ell) \eqs,
\eeq
for all $k \geq 0$. Observe that $\bar{R}(k) \geq R(k)$ for all $k \geq 0$
and it follows from the first part that $\bar{R}(k) \leq \delta^k \bar{R}(0)$. 
Now,
\beq \notag
	\begin{split}
		& \frac{\bar{R}(k+1)}{\bar{R}(k)} = p + \frac{q m_1^{(k)} \max_{ (k-m_2^{(k)})_+ \leq \ell \leq k } \bar{R}^\eta ( \ell )}{ \bar{R}(k) } \\
		& \le p + q m_1^{(k)}  \max_{ (k-m_2^{(k)})_+ \leq \ell \leq k }\left\{\frac{ \bar{R} ( \ell )}{\bar{R}(k)}\delta^{(\eta-1)\ell}\right\} \bar{R}(0)^{\eta-1}\\
		& \le p + q m_1^{(k)}  \max_{ k-m_2^{(k)} \leq \ell \leq k }\left\{\frac{R (\ell )}{R(k)}\right\}\delta^{(\eta-1)(k-m_2^{(k)})}R(0)^{\eta-1}\\	
	\end{split}
\eeq
Using the fact that $\bar{R}(k+1) / \bar{R}(k) \geq p$ shows
\beq \notag
\frac{\bar{R}(k+1)}{\bar{R}(k)} \leq p + q R(0)^{\eta-1} m_1^{(k)} p^{-m_2^{(k)}} \delta^{(\eta-1)(k-m_2^{(k)})} 
\eeq
Letting $\eta' = \eta - 1$, we see that
\beq
\begin{split}
& m_1^{(k)} p^{-m_2^{(k)}} \delta^{\eta'(k-m_2^{(k)})} = m_1^{(k)} (p \delta^{\eta'})^{-m_2^{(k)}} \delta^{\eta' k }\\
& < m_1^{(k)} (p \delta^{\eta'})^{k-m_2^{(k)}} \to 0,
\end{split}
\eeq
as $k\to \infty$, under the assumptions on $m_1^{(k)}, m_2^{(k)}$.  
Finally, taking the limit $k \rightarrow \infty$ shows:
\beq
\lim_{k \rightarrow \infty} \bar{R}(k+1) / \bar{R}(k) = p \eqs.
\eeq

\ifplainver
\section{Feasible Regions of Step Size} \label{app:step}
\else
\subsection{Feasible Regions of Step Size} \label{app:step}
\fi
Let $\eta \eqdef \min_j \eta_j$. 
Under A\ref{ass:delay}, the required condition \eqref{eq:delay_gen} in Lemma~\ref{lem:nonlinear} 
is satisfied. 
Moreover, a lower bound to $\xi_j^\star(\delta)$ for all $j$ can be evaluated --- 
following \eqref{eq:xi_lb} and solving 
the convex minimization in the bound, it can be shown 
that
\beq 
\xi^\star(\delta) \geq - \frac{ 1 }{\log \delta} + \frac{ c_0 + \log( 1 / (\beta \log(1/\delta)) ) }{\log \delta} {-\eta m_0}
\eeq
Consequently, using  the inequality $\log(x) \geq 1 - 1/x$ gives:
\beq
\begin{split}
\delta^{-\xi^\star(\delta)} & \geq e^{ 1 - c_0 - \log( 1 / (\beta \log(1/\delta))  ) }\delta ^{ \eta m_0} \\
& = e^{ 1 - c_0} \delta^{\eta m_0} {\beta} \log(1/\delta)\\
& \geq e^{ 1 - c_0} \delta^{\eta m_0} \beta  (1 - \delta ) 
\end{split}
\eeq
Taking $\delta = 1 - \gamma \Delta \frac{2 \mu L}{\mu + L}$, the sufficient condition 
\eqref{eq:lem_1_pre} is satisfied
given the following conditions hold:
\beq
\gamma^{ 3 \lceil \frac{j}{2} \rceil } \leq C_j' 
\delta^{\eta m_0} \cdot \gamma^2,~j=1,...,4 \eqs.
\eeq
Now let us define $A \eqdef \frac{2 \mu L}{\mu + L}$. 
If we consider the case when $j=1,2$, we have
\beq
\gamma \leq C_j' 
(1 - \Delta A {\gamma} )^{\eta m_0} 
\eeq
Applying Bernoulli's inequality 
gives the following \emph{sufficient} condition:
\beq
\gamma \leq C_j'  ( 1 - \eta m_0 \Delta A  {\gamma} )  
\eeq
Therefore, we have
\beq \label{eq:bdbd}
\begin{split}
\gamma & \leq \frac{ C_j'  }{ 1 +  C_j' \Delta \eta m_0 A }  = \Big( \Delta \eta m_0 A + \frac{1}{ C_j'} \Big)^{-1},~j=1,2.
\end{split}
\eeq
Repeating the same analysis for $j=3,4$ yields
the sufficient condition
\beq
\gamma^3 \leq C_j' 
(1 -\eta m_0  \Delta A {\gamma} ) 
\eeq
Using the upper bound $\gamma^2 \leq 1$ gives the same conditions
as in \eqref{eq:bdbd}:
\beq \label{eq:lasteq}
\gamma \leq \Big( \Delta \eta m_0 A + \frac{1}{C_j'} \Big)^{-1},~j=3,4.
\eeq

\fi

\end{document}

%% file: Tikz/SUCAG.tikz.tex
\begin{tikzpicture}[scale=0.7,thick]
  \foreach \place/\name in {{(0,-1)/a}, {(2.5,-.5)/b}, {(1.25,3.25)/f}, {(-3.5,2)/g}}
    \node[superpeers, black, fill=white, opacity=.3] (\name) at \place {};
  \node (n1) at (-5,0) {};
  \node (n2) at (-4,-1.5) {};
  \node (n3) at (4,3) {};
  \node (n4) at (4,0) {};
  \node (n5) at (-5,2.5) {}; 
  \node (n6) at (3.25,3.5) {};
  \foreach \place/\name in {{(-3,-0.5)/e}, {(-1,2)/d}, {(2.5,1.5)/c}}
    \node[superpeers, asublue!70!black, very thick, fill=white] (\name) at \place {};
  \foreach \source/\dest in {a/b, b/a, a/c, c/a, a/d, b/c,a/e, e/a, d/f, f/c, e/g, g/d, g/f, e/n1, e/n2, c/n3, c/n4, g/n5, f/n6}
    \path[color=gray!40] (\source) edge (\dest);
  \draw[very thick, asublue!70!black,->] (e) to node[auto, yshift=-.2cm, xshift=.1cm] {\scriptsize ${\bm H}^{k-1}, {\bm b}^{k-1}, \prm^k$}  (d);
  \draw[very thick, asublue!70!black,->] (d) to node[auto, xshift=-.2cm] {\scriptsize ${\bm H}^{k}, {\bm b}^{k}, \prm^{k+1}$} (c);
  \node [below, yshift=-0.1cm, xshift=.3cm] at (e) {\tikz{\node[fill=white] {\scriptsize Agent $i_{k-1}$}}};
  \node [below, yshift=-0.1cm, xshift=.65cm] at (d) {\tikz{\node[fill=white] {\scriptsize Agent $i_{k}$}}};
  \node [below, yshift=-0.1cm, xshift=.4cm] at (c) {\tikz{\node[fill=white] {\scriptsize Agent $i_{k+1}$}}};
  \node [text width = 3.75cm,above, yshift=.6cm, xshift=-.8cm, legend_graph] (cap) at (d) {\scriptsize 1.~Compute \eqref{eq:sucag}, \eqref{eq:sucag_a}, \eqref{eq:inc_upd}.\\[-.05cm] 2.~Store $\prm^k$ in local memory.};
  \draw[->, dashed, very thick] (cap) -- (d);
\end{tikzpicture}

%% file: Tikz/cdc_test_case.tikz.tex
\resizebox{.98\linewidth}{!}
{ \sf 
\begin{tikzpicture}
\pgfplotsset{ scale only axis,
    width=0.55\textwidth,height=0.4\textwidth,
    grid=both,grid style={line width=.01pt, draw=gray!30},
    legend cell align=left, legend style={legend pos=north east,font=\large, yshift = -0.0cm, opacity  = 0.8},
    xlabel={\Large Iteration number $\times N$},
    enlarge x limits=0.05,
    enlarge y limits=0.05,
    xmin = 1.1, xmax = 95, ymin = 1e-12
}

\pgfplotsset{every tick label/.append style={font=\large}}

\begin{semilogyaxis}[ylabel shift = -0.1cm, ylabel={\Large Opt.~Gap $F(\prm^k) - F(\prm^\star)$}]
\addplot[blue, very thick, mark options={solid,mark size=4}, mark=diamond, each nth point={2}] 
      table[x index=0, y index=1, col sep=comma] {./Tikz/CDC_step_01L_curv.csv};
      \addlegendentry{Determ.~CIAG ($\gamma=L^{-1}$)};
\addplot[asublue, very thick, mark options={solid,mark size=4}, mark=triangle, each nth point={5}] 
      table[x index=0, y index=1, col sep=comma] {./Tikz/CDC_step_01L_sag.csv};
      \addlegendentry{D-SAG ($\gamma = 0.2/L$)};
\addplot[asuorange!80!black, very thick, mark options={solid,mark size=3}, mark=o, each nth point={5}] 
      table[x index=0, y index=1, col sep=comma] {./Tikz/CDC_step_01L_dgd.csv};
      \addlegendentry{DIG-ing ($\gamma = 0.02/L$)};
\addplot[asured, very thick, mark options={solid,mark size=3}, mark=square, each nth point={5}] 
      table[x index=0, y index=1, col sep=comma] {./Tikz/CDC_step_01L_curv_rand.csv};
      \addlegendentry{D-CIAG ($\gamma = 0.1/L$)};
\addplot[black, very thick, mark options={solid,mark size=3}, mark=*, each nth point={5}] 
      table[x index=0, y index=1, col sep=comma] {./Tikz/CDC_step_01L_ucag_rand.csv};
      \addlegendentry{D-SUCAG ($\gamma = 0.1/L$)};
\addplot[asured, dashed, very thick, mark options={solid,mark size=3}, mark=square, each nth point={5}] 
      table[x index=0, y index=1, col sep=comma] {./Tikz/CDC_step_04L_curv_rand.csv};
      \addlegendentry{D-CIAG ($\gamma = 0.4/L$)};
\addplot[black, dashed, very thick, mark options={solid,mark size=3}, mark=*, each nth point={5}] 
      table[x index=0, y index=1, col sep=comma] {./Tikz/CDC_step_04L_ucag_rand.csv};
      \addlegendentry{D-SUCAG ($\gamma = 0.4/L$)};
\end{semilogyaxis}

\end{tikzpicture}
} 